\documentstyle[12pt]{article}
\newtheorem{theorem}{{\sc Theorem}}
\newcommand{\bt}{\begin{theorem}}
\newcommand{\et}{\end{theorem}}
\newcommand{\newsection}[1]{\setcounter{equation}{0} \setcounter{theorem}{0}
\section{#1}}

\newcommand{\NI}{\noindent}

\newcommand{\bea}{\begin{eqnarray}}
\newcommand{\eea}{\end{eqnarray}}

\def \b #1 {\bf #1}

\newcommand{\cla}{{\cal A}}

\newcommand{\cls}{{\cal S}}

\newcommand{\clh}{{\cal H}}

\newcommand{\clb}{{\cal B}}

\newcommand{\cll}{{\cal L}}

\newcommand{\raro}{\rightarrow}

\newcommand{\vsp}{\vskip 1em}

\def \qed {\hfill \vrule height6pt width 6pt depth 0pt}
\newcommand{\be}{\begin{equation}}
\newcommand{\ee}{\end{equation}}
\newcommand{\ben}{\begin{eqnarray*}}
\newcommand{\een}{\end{eqnarray*}}
\begin{document}
\thispagestyle {empty}
\sloppy

\centerline{\large \bf A resolution of quantum dynamical semigroups }

\bigskip
\bigskip
\centerline{\bf Anilesh Mohari }
\smallskip
\centerline{\bf S.N.Bose Center for Basic Sciences, }
\centerline{\bf JD Block, Sector-3, Calcutta-91 }
\centerline{\bf E-mail:anilesh@boson.bose.res.in}
\bigskip
\smallskip
\centerline{\bf Abstract}
\smallskip
\bigskip
\vsp
We consider a class of quantum dissipative systems governed by a 
one parameter completely positive maps on a 
von-Neumann algebra. We introduce a notion of recurrent and 
metastable projections for the dynamics and prove that the unit
operator can be decomposed into orthogonal projections where
each projections are recurrent or metastable for the dynamics.    

\hfil\eject

\newsection{ Introduction:}

\vsp
Let $\tau=(\tau_t,\;t \ge 0)$ be a semigroup of identity preserving completely 
positive maps [Da1,Da2,BR] on a von-Neumann algebra $\cla_0$ acting on a 
separable Hilbert space $\clh_0$, where either the parameter $t \in \!R_+$, the 
set of positive real numbers or $ t \in \!Z^+$, the set of positive integers. 
In case $\cla_0 =l^{\infty}(\cls)$, the bounded measurable functions on a 
countable set $\cls$, $\tau$ is the Markov semigroup associated with a 
continuous or discrete time Markov chain. In such a case the state space 
$\cls$ is well classified as disjoint union of recurrent and metastable 
states. Such a classification enriched our understanding of asymptotic 
behavior of associated stochastic processes. In this exposition we aim 
to achieve a similar classification in the general framework of 
non-commutative probability and prove results related to the 
asymptotic behavior of the quantum dynamical systems $(\tau_t)$ 
if it admits a normal invariant state. To that end
we assume further that the map $\tau_t$ is normal for each $t \ge 0$ 
and the map $t \raro \tau_t(x)$ is weak$^*$ continuous for each $x \in \cla_0$.     

\vsp
Following [FR3] we say a projection $p \in \cla_0$ is {\it sub-harmonic and harmonic  } if $\tau_t(p) \ge p$ and $\tau_t(p)=p$ 
for all $t \ge 0$ respectively. For a sub-harmonic projection $p$, we 
define the reduced quantum dynamical semigroup $(\tau^p_t)$ on the von-Neumann 
algebra $p\cla_0p$ by $\tau^p_t(x)=p\tau_t(x)p$ where $t \ge 0 $ and $x \in \cla^p_0$. 
If $p$ is also a harmonic projection, the reduced dynamics $\tau_t^{1-p}$ and
$\tau^p_t$ are complement to each other. For more details we refer to [Ev]. Following Evans [Ev] 
we also say $(\tau_t)$ is {\it irreducible } if there exists no 
non-trivial harmonic projection in $\cla_0$ for $(\tau_t)$. This notion is 
not quite equivalent to the notion of irredcibility for Markov chain [No].   

\vsp
A non-zero sub-harmonic projection $p$ is said to be {\it minimal} if there exists no 
sub-harmonic projection $q$ so that $0 < q < p$. A non-zero projection 
$p \in \cla$ is called {\it recurrent } if $p$ is minimal  
sub-harmonic. A recurrent projection $p$ is called {\it positive recurrent} 
if there exists a normal invariant state for $(\tau_t)$ with support equal 
to $p$. A recurrent projection is called {\it null
recurrent } if it is not a positive recurrent. If the identity operator 
is a recurrent projection, we say $(\tau_t)$ is recurrent.   

\vsp
For a sub-harmonic projection $p$, $1$ is an upper bound for the increasing positive 
operators $\tau_t(p),\; t \ge 0$. Thus there exists an operator $0 \le y \le 1$ so that 
$y=\mbox{s.lim}_{t \raro \infty} \tau_t(p)$. If $y$ is an injective map on 
the Hilbert space $\clh_0$ we say the projection $1-p$ is {\it metastable }.  
If $y=1$ we say $1-p$ is {\it transient}. A normal state $\phi_0$ is called 
{\it invariant} for $(\tau_t)$ if $\phi_0 \tau_t(x) = \phi_0(x)$ for all $x \in 
\cla_0$ and $t \ge 0$. The support $p$ of a normal invariant state is a sub-harmonic 
projection [FR3] and $\phi^p_0$, the restriction of $\phi_0$ to $\cla^p_0$ is 
an invariant normal state for $(\tau^p_t)$. Thus asymptotic properties ( ergodic, mixing ) 
of the dynamics $(\cla_0,\tau_t,\phi_0)$ is well determined by the asymptotic properties 
(ergodic, mixing respectively ) of the reduced dynamics $(\cla^p_0,\tau^p_t,\phi^p_0)$ 
provided $y=1$. For more details we refer to [Fr,Mo1].  

\vsp
We explore a necessary and sufficient condition for a sub-harmonic
projection to be metastable and also find a sufficient condition for
a metastable projection to be transient. In case the von-Neumann
algebra is type-I with centre completely atomic, we will show that
there exists a mutually orthogonal family of recurrent projections 
$(p_i)$ so that $\sum_i p_i \le 1$ and $1-\sum_ip_i$ is a 
metastable projection. However such a decomposition is not 
unique in general unlike in the classical Markov chain. Thus the problem 
remain open also for more general non-commutative algebra of observables.  
Related results for stochastically complete Brownian flows on a smooth 
manifold associated with a family of vector fields are proved in [Mo2].

\vsp
In particular for any finite dimensional Hilbert space $\clh_0$ we prove 
that any metastable projection is also transient and any recurrent projection 
is also positive recurrent. Thus even in these general non-commutative 
dynamical systems of irreversible processes classical feature remain valid 
once the Hilbert space is finite dimensional. If there exists a 
unique recurrent projection $p$ and $\clh_0$ is finite dimensional then
there exists an unique invariant normal state and $(\tau_t)$ is ergodic 
with respect to that state.    
 
\vsp
\newsection{ Sub-harmonic projections: }

\bigskip
In the following we list few crucial property ( Proposition 2.1 and Proposition 2.2 ) of sub-harmonic projections 
appeared in [FR3,Mo1].

\vsp
\NI {\bf PROPOSITION 2.1: } Let $p$ be a sub-harmonic projection 
for $(\tau_t)$. Then the following hold:

\NI (a) for all $t \ge 0$, $p\tau_t(p) = \tau_t(p)p=p$.

\NI (b) $\tau_t(x(1-p))p=0$ for all $x \in \cla_0, t \ge 0$. 

\vsp
For a projection $p$, $\cla^p_0=p\cla_0p$ is also a von-Neumann algebra 
acting on the Hilbert subspace $p\clh_0$. Thus for a sub-harmonic projection 
$p$ we verify by Proposition 2.1 that $(\tau^p_t)$ defined by 
$\tau^p_t(x)= p\tau_t(x)p,\; x \in \cla^p_0$ is a Markov semigroup. 
Let the strong limit of $\tau_t(p)$ be $y$ as $t \uparrow \infty$. By 
Proposition 2.1 (a) we have $py=yp=p, p \le y \le 1$ and $\tau_t(y)=y \forall 
t \ge 0$. Thus $p\tau_t(1-y^2)p=p\tau_t(p(1-y^2)p)p=0$. So we also have 
$p\tau_t(y^2) = \tau_t(y^2)p=p$ for all $t \ge 0$. Since $\tau_t(p) \le 
\tau_t(y^2) \le \tau_t(y)=y$, the strong limit of $\tau_t(y^2)$ as $t \raro \infty$ is also 
$y$. In case $y^2$ is also an invariant element for $(\tau_t)$, we have 
$y^2=y$. In general $y^2$ need not be an invariant element even for an 
irreducible classical Markov semigroup $(\tau_t)$. For more details we refer to [Mo1]. 

\vsp
\NI {\bf PROPOSITION 2.2: } Let $p$ be a sub-normal projection and 
$y=\mbox{s.lim}\tau_t(p)$. Then for any $z \in \clb(\clh_0)$ the following 
statements are equivalent:

\NI (a) $yz=0$

\NI (b) $\tau_t(p)z=0$ for all $t \ge 0$.
                               
\bigskip
The following proposition is crucial for classification and related 
problems.

\vsp
\NI {\bf PROPOSITION 2.3:} Let $x \in \cla_0$ be a non-negative invariant 
element for a completely positive unital normal map $\tau$ on a von-Neumann algebra 
$\cla_0$ and $q$ be the projection on the closure of the range of $x$. Then $1-q$ is a 
sub-harmonic projection for $\tau$.  

\NI {\bf PROOF:} Let $q$ be the unique minimal projection so that $(1-q)x=0$.
For any unitary element $u' \in \cla'_0$ in the commutant, we have 
$(1-q)u'xu'^*=0$, thus $u'^*(1-q)u'x=0$. Since $u'^*(1-q)u'$ is also a 
projection, by minimality we conclude that $u'^*qu' \ge q$. $u'$ being an 
arbitrary unitary element in $\cla'_0$, we conclude that $u'^*qu'=q$. That 
$q$ is an element in $\cla_0$ follows now by von-Neumann's density theorem 
[BR]. By Steinspring representation [BR], there exists a 
countable family of bounded operators $(l_k)$ so that $\tau(z) = \sum_k l_k^*z
l_k$ for all $z \in \cla_0$. Since 
$x$ is an invariant element, we have $(1-q)\tau(x)=(1-q)x=0$. Hence we have 
by the representation $xl_k(1-q)=0$. Thus we have $ql_k(1-q)=0$. Once more by 
the representation we conclude that $\tau(q)(1-q)=0$, i.e. $\tau(q) \le q$. \qed 

\vsp
For the classification we need to guarantee existence of a minimal projection. 

\vsp
\NI {\bf THEOREM 2.4:} Let $\cla_0$ be type-I with centre 
completely atomic. There exists a commuting family $(p_i)$ of 
countable many orthogonal recurrent projections for $(\tau_t)$ so 
that $q=1-\sum_ip_i$ is metastable.

\NI {\bf PROOF:} Since $\cla_0$ is type-I with centre completely atomic,
minimal projection exits [BR,vol-1]. In case there exists no non-trivial 
minimal sub-harmonic projection, we have trivial partition $p_1=1$ and $q=0$. 
Suppose $p_1$ is a minimal sub-harmonic projection and $q_1$ is the 
projection on the closure of the range of $y_1$ where $y_1$ is the strong limit of
$\tau_t(p_1)$ as before. Thus by Proposition 2.3 $1-q_1$ is sub-harmonic.
In case $1-q_1=0$ or $1-q_1$ is recurrent then $1-p_1-p_2$ is metastable,
where $p_2=1-q_1$. So the construction ends. Otherwise we look for a minimal
sub-harmonic projection $0 < p_2 < (1-q_1)$ and so on. Note that by
our construction the family is countable ($\clh_0$ is separable) and 
commuting.  \qed

\bigskip
\NI {\bf THEOREM 2.5:} Let $p$ be a sub-harmonic projections for $(\tau_t)$ so that 
$y \ge \lambda$ for some $\lambda > 0$. Then $y=1$.
 
\smallskip
\NI {\bf PROOF:} Since $yp=p$, by our hypothesis we write
$y \ge p + \lambda (1-p)$ for some $\lambda > 0$. 
Thus for any $t \ge t_0$ we have $y \ge \tau_t(p) + \lambda (1-\tau_t(p))$. 
By taking limit $t \raro \infty$ we conclude that 
$y \ge y + \lambda(1- y)$. Hence $y \ge 1$, since $\lambda > 0$. 
Since $0 \le y \le 1$ we conclude that $y=1$. \qed

\bigskip
\NI {\bf COROLLARY 2.6:} Let $p$ be a maximal non-trivial sub-harmonic
projection $p \neq 1$ for $(\tau_t)$ and $1-p$ be finite dimensional. 
Then $\tau_t(p) \raro 1$ as $t \raro \infty$ in strong operator topology.
 
\smallskip
\NI {\bf PROOF:} Let $q$ be the projection on the range of $y$. Since $y$ is
a non-negative invariant element for $(\tau_t)$, Proposition 2.3 says that
$1-q$ is sub-harmonic. Since $p$ is a maximal non-trivial sub-harmonic projection
which commutes with $1-q$, we have $1-q=0$. Hence $y$ is one to
one. Since $1-p$ is finite dimensional, $y$ is bounded away from
zero i.e. there exists a constant $\lambda > 0$ so that $y \ge \lambda 1$.
Thus the proof is complete once we appeal to Theorem 2.5. \qed

\vsp 
If the parameter $t \in \!Z_+$, the dynamics is determined 
by a unital completely positive normal map $(\tau)$ and $\tau_n=
\tau \circ ..\circ \tau$- $n$ fold composition. In such a case $p$
is sub-harmonic for $(\tau_n)$ if and only if $\tau(p) \ge p$. 

\vsp
\NI {\bf THEOREM 2.7:} Let $\tau$ be a unital completely positive 
map on a von-Neumann algebra $\cla_0$. Then $p$ be a sub-harmonic 
projection for $\tau$ if and only if $(1-p)l_kp=0$ for all 
$k \ge 1$, where $(l_k)$ is the family of bounded operators 
so that $\tau(x)=\sum_k l_k^*xl_k$ for all $x \in \cla_0$.    

\NI Moreover for a sub-harmonic projection $p$ for $\tau$, 
and a bounded operator $z$, $yz=0$ if and only if $pl_{i_1}l_{i_2}l_{i_3}
...l_{i_n}z=0$, where $s.\mbox{lim}_{n \raro \infty}\tau_n(p)=y$. 

\vsp
\NI {\bf PROOF:} Let $p$ be sub-harmonic for $\tau$, i.e.
$\tau(p) \ge p$. We have $p\tau(p)p \ge p$, $p$ being 
a projection. Hence $p\tau(1-p)p \le 0$. Since $1-p \ge 0$,
we also have $p\tau(1-p)p \ge 0$. So $p\tau(1-p)p=0$. By
the representation we have $\sum_k pl^*_k(1-p)l_kp=0$, 
hence the result follows. The last part follows once 
we note that $yz=0$ if and only if $\tau_n(p)z=0$ for all $n \ge 0$. \qed 

\vsp
Hence $y$ is injective if and only if $p$ together with 
the range of these $l^*_{i_1}l^*_{i_2}..l^*_{i_m}p$ family 
of operators, where $i_1,i_2,..i_m \ge 1$ generate the Hilbert space $\clh_0$. 
In such a case for finite dimensional Hilbert space $\clh_0$, Theorem 2.5 
says that $y=1$. Although the family $(l_k)$ is not uniquely determined
by the map $\tau$, the criteria appeared in Theorem 2.7 is independent of
the choice we make for it's representation.

\newsection{ Strong ergodic property: }

\vsp
In the last section we investigated various criteria for the support projection $p$ associated with 
the normal invariant state so that $1-p$ is transient for the dynamics $(\cla_0,\tau_t,\phi_0)$. In a 
recent paper [Mo1], we have investigated how various asymptotic properties ( ergodic, weak mixing, strong mixing ) 
of $(\cla_0,\tau_t,\phi_0)$ can be determined by the reduced dynamics $(\cla_0^p,\tau^p_t,\phi_0)$ when $1-p$ is transient. 
In the following theorem we prove one more such a property known in the literature as strong ergodicity.

\vsp
\NI {\bf THEOREM 3.1:} Let $(\cla_0,\tau_t,\phi_0)$ be a quantum dynamical systems with a 
normal invariant state $\phi_0$ with support projection $p$. If $\mbox{s-limit}_{t \raro \infty}\tau_t(p)=1$
then the following statements are equivalent:

\NI (a) $||\phi\tau_t-\phi_0|| \raro 0$ as $t \raro \infty$ for any normal state on $\cla_0$.

\NI (b) $||\phi \tau^p_t- \phi_0|| \raro 0$ as $t \raro \infty$ for any normal state on $\cla_0$.

\vsp
\NI {\bf PROOF: } That (b) implies (a) is trivial. For the converse we write 
$||\phi \tau_t-\phi_0|| = \mbox{sup}_{x:||x|| \le 1} |\phi \tau_t(x)-\phi_0(x)| 
\le \mbox{sup}_{\{x:||x|| \le 1\}} |\phi \tau_t(pxp)-\phi_0(pxp)| +
\mbox{sup}_{ \{x:||x|| \le 1 \} } |\phi \tau_t(pxp^{\perp})| +
\mbox{sup}_{ \{x:||x|| \le 1 \} } |\phi \tau_t(p^{\perp}xp)| +
\mbox{sup}_{ \{x:||x|| \le 1 \} } |\phi \tau_t(p^{\perp}xp^{\perp})|$.
Since $\tau_t((1-p)x) \raro 0$ in the weak$^*$ topology and 
$|\phi \tau_t(xp^{\perp})|^2 \le |\phi \tau_t(xx^*)|\phi(\tau_t(p^{\perp}))| \le ||x||^2 \phi(\tau_t(p^{\perp})$
it is good enough if we verify that (a) is
equivalent to  $\mbox{sup}_{\{ x:||x|| \le 1 \} } |\phi\tau_t(pxp)-\phi_0(pxp)| \raro 0$
as $t \raro \infty$. To that end we first note that limsup$_{t \raro
\infty} \mbox{sup}_{x:||x|| \le 1} |\psi(\tau_{s+t}(pxp))-\phi_0(pxp)|$ is independent 
of $s \ge 0$ we choose. On the
other hand we write $\tau_{s+t}(pxp) = \tau_s(p\tau_t(pxp)p) +
\tau_s(p\tau_t(pxp)p^{\perp})+ \tau_s(p^{\perp}\tau_t(pxp)p) +
\tau_s(p^{\perp}\tau_t(pxp)p^{\perp})$ and use the fact for any
normal state $\phi$ we have
$\mbox{limsup}_{t \raro \infty}\mbox{sup}_{x:||x|| \le 1}
|\psi(\tau_s(z\tau_t(pxp)p^{\perp})|\le ||z||\;|\psi(\tau_s(p^{\perp}))| $ for all $z \in \cla_0$.
Thus by our hypothesis on the support and we conclude that (a) hold whenever (b) is true. \qed

\vsp
\newsection{ Quantum mechanical master equation:}

\vsp
In this section we aim to deal with a class of quantum dynamical semigroup.
We say a normal Markov semigroup $(\tau_t)$ on $\cla_0$ is norm continuous if
limit$_{t \raro 0 }||\tau_t-I||=0$. In such a case the generator $\cll$
is a bounded operator on $\cla$ and can be described [GKS,Li,CE] by
\be
\cll(x)= Y^*x +xY + \sum_{k \ge 1} L_k^*xL_k
\ee
where $Y \in \cla_0$ is the generator of a norm continuous
contractive semigroup on $\clh_0$ and $L_k,\;k \ge 1$ is a
family of bounded operators so that $\sum_k L^*_kxL_k \in
\cla_0$ whenever $x \in \cla_0$. However this choice $(Y,L_k,\;k
\ge 1)$ is not unique. Conversely, for any such a family $(Y,L_k)$
with $Y \in \cla_0$ and $L_kxL^*_k \in \cla_0,\;\forall x \in \cla_0$, there
exists a unique Markov semigroup $(\tau_t)$ with $\cll$ as its generator.
There are many methods to show the existence of a Markov 
semigroup $(\tau_t)$ with $\cll$ as it's generator [Da3,MS,CF].
Here we describe one such a method [CF].

\vsp
We consider the following iterated equation :

\vsp
\begin{eqnarray}
\tau^0_t(x) & = & e^{tY^*}xe^{tY}  \\
\tau^{(n)}_t(x) & = & e^{tY^*}xe^{tY} + \int^t_0e^{(t-s)Y^*}\Phi(
\tau^{(n-1)}_s(x))e^{(t-s)Y}ds
,\;n \ge 1 
\end{eqnarray}
where $\Phi(x)=\sum_k L^*_kxL_k$. It is simple to check for $x \ge 0$ that
$$ 0 \le \tau^{n-1}_t(x) \le \tau^n_t(x) \le ||x||\;I,\;\forall t \ge 0$$
Thus we set for $x \ge 0,\;\tau_t(x)=limit_{n \raro \infty} \tau^{(n)}_t(x)$
in the weak$^*$ topology. For an arbitrary element we extend it by linearity. 
Thus we have 
\be
\tau_t(x) =  e^{tY^*}xe^{tY} + \int^t_0e^{(t-s)Y^*}\Phi(\tau_s(x))e^{(t-s)Y}ds
\ee
for any $x \in \cla_0$.

\vsp
In such a case [Ev], it is simple to check that $(\tau_t)$ is irreducible
if and only if $\{x \in \cla_0:[x,H]=0,\;[x,L_k] = 0 \forall k \ge 1 \}$ 
is trivial. 

The following simple but important result due to Fagnola-Rebolledo [FR3].

\vsp
\NI {\bf THEOREM 4.1:} A projection $p$ is sub-harmonic if and only if 
$(1-p)Yp=0$ and $(1-p)L_kp=0$ for all $1 \le k \le \infty$. 

\vsp
\NI {\bf PROOF:} For a proof and a more general result we 
refer to [FR3]. \qed

\vsp
\NI {\bf PROPOSITION 4.2:} Let $p$ be a sub-harmonic projection and $y=\mbox{s.lim}_
{t \raro \infty}\tau_t(p)$. For any $z \in {\cal B}(\clh_0)$ following are 
equivalent:

\NI (a) $yz=0$ 

\NI (b) $pz=0,\;pL_{i_1}L_{i_2}....L_{i_n}z=0$ for all $0 \le i_m < \infty$, 
$1 \le m \le n$ and $n \ge 1$, where $L_0=Y$. 

\vsp
\NI {\bf PROOF:} By Proposition 2.2, $yz=0$ if and only if $z^*\tau_t(p)z=0$ for all 
$t \ge 0$. Now by (4.4) we have $z^*\tau_t(p)z=0$ if and only if 
$z^*e^{tY^*}pe^{tY}z=0$ and 
$z^*\Phi(\tau_t(p))z=0$ for $t \ge 0$. Thus we have $pe^{tY}z=0$ and also 
$z^*\Phi(\tau_t(p))z=0$ for all $t$. Thus in particular we have 
$z^*\Phi(p)z=0$, hence $pL_kz=0$ for all $k \ge 0$. We go now by induction 
on $n$, we check if $z'=pL_{i_1}L_{i_2}....L_{i_n}z$ then $yz'=z'$ thus (a) 
implies (b). For the converse statement, we check that derivative of any 
order at $t=0$ of $z^*\tau_t(p)z$ vanishes, thus constant which is zero. \qed   

\vsp
Thus the zero operator is the only element $z$ that satisfies (b) 
if and only if the closure of the range of $y$ is the entire Hilbert 
space. Thus $p$ together with $L^*_{i_1}L^*_{i_2}....L^*_{i_n}p$ where 
$0 \le i_m < \infty$, $1 \le m \le n$ and $n \ge 1$ will generate the 
Hilbert space if and only if $y$ is one to one. In particular we find this 
property is a necessary condition for $y$ to be $1$. In general the 
condition is not a sufficient one. Once more one can construct a counter 
example in birth and death processes where no population is 
an absorbing state and birth and death rates are such that the population will 
extinct with positive probability but need not be $1$. We omit the details. 
However Corollary 2.6 says that for finite dimensional $1-p$  
the condition is also sufficient for $y=1$. So we have the following result. 

\vsp
\NI {\bf THEOREM 4.3:} Let $p$ be a sub-harmonic projection for $(\tau_t)$ 
and $1-p$ be finite dimensional. Then $y=1$ if and only if 
$p$ together with the ranges of operators $L^*_{i_1}....L^*_{i_n}p$,
$i_i,..,i_n \ge 0$ generates the entire Hilbert space $\clh_0$.   

\vsp
Although operators $(L_k)$ are not uniquely determined for a given 
$(\tau_t)$, the criteria appeared in Theorem 4.3 is independent of choice
we make. We end this exposition with a comment that the result similar to 
Proposition 4.2 with appropriate domain condition [FR3] can be proved
to include a class of quantum irreversible system where the $\cla_0={\cal B}
(\clh_0)$ with generator which need not be a bounded operator.

\bigskip
{\centerline {\bf REFERENCES}}

\begin{itemize} 

\bigskip
 
\item {[BR]} Bratelli, O., Robinson, D.W. : Operator algebras
and quantum statistical mechanics, I,II, Springer 1981.

\item {[CF]} Chebotarev, A.M., Fagnola, F. Sufficient conditions for
conservativity of minimal quantum dynamical semigroups, J. Funct. Anal.
153 (1998) no-2, 382-404.

\item{[CE]} Christensen, E., Evans, D. E.: Cohomology of
operator algebras and quantum dynamical semigroups, J.Lon. Maths
Soc. 20 (1970) 358-368.

\item{[Da1]} Davies, E.B.: One parameter semigroups, Academic
Press, 1980.

\item{[Da2]} Davies, E.B.: Quantum theory of open systems,
Academic press, 1976.

\item{ [Da3}] Davies, E.B.: Quantum dynamical semigroups and the neutron diffusion equation. Rep. Math. Phys. vol-11, no-2 (1977), 169-188. 

\item{[Ev]} Evans, D.E.: Irreducible quantum dynamical
semigroups, Commun. Math. Phys. vol-54 (1977), 293-297.

\item{[FR1]} F.~Fagnola, R.~Rebolledo.
The approach to equilibrium of a class of quantum dynamical semigroups
Inf. Dim. Anal. Q. Prob. and Rel. Topics, vol-1, no. 4 (1998),  1--12.
 
\item {[FR2]} F.~Fagnola and R.~Rebolledo.
On the existence of invariant states for quantum dynamical semigroups
J. Math. Phys., vol-42 (2001), 296-1308.
 
\item {[FR3]} F.~Fagnola, R.~Rebolledo.: Subharmonic projections for a Quantum Markov Semigroup, 
J. Math. Phys., vol-43, no. 2 (2002), 1074--1082.
 
\item{[Fr]} Frigerio, A.: Stationary states of quantum dynamical
semigroups, Commun. Math. Phys. vol-63 (1978), 269-276.

\item{[GKS]} Gorini, V., Kossakowski, A., Sudarshan, E.C.G. : Completely
positive dynamical semigroups of n-level systems, J. Math. Phys. vol-17 
(1976) 821-825.

\item{[Li]} Lindblad, G. : On the generators of quantum dynamical semigroups, 
Commun.  Math. Phys. vol-48 (1976), 119-130.

\item{[Mo1]} Mohari, A.: Markov shift in Non-commutative probability, J. of Funct. Anal. vol-199, 
no. 1 (2003), 189--209.  

\item{[Mo2]} Mohari, A.: Ergodicity of Homogeneous Brownian flows, Stochastic processes and their 
application vol-105, no. 1 (2003), 99--116.

\item{[MS]} Mohari, A., Sinha, K.B.: Stochastic dilation of minimal quantum 
dynamical semigroup. Proc. Indian Acad. Sci. Math Sci. vol-102 (1992), no-2, 
159-173. 

\item{[No]} Norris, J.R., Markov Chain, Cambridge University Press, 1998.

\end{itemize}

\end{document}